\documentclass[graybox]{svmult}

\usepackage{type1cm}        % activate if the above 3 fonts are
\usepackage{makeidx}         % allows index generation
\usepackage{graphicx}        % standard LaTeX graphics tool
\usepackage{multicol}        % used for the two-column index
\usepackage[bottom]{footmisc}% places footnotes at page bottom

\usepackage{amsmath,amssymb}
\usepackage{booktabs}
\providecommand{\mat}[1]{\ensuremath { {\rm \bf #1} } }
\providecommand{\restrict}[2]{\ensuremath \left. #1 \right|_{#2}}
\providecommand{\surf}{\ensuremath \mathcal{S} }
\providecommand{\pdiff}[2]{\ensuremath \frac{\partial #1}{\partial #2} }
\DeclareMathOperator*{\argmin}{arg\,min}
\usepackage{newtxtext}       %
\usepackage{newtxmath}       % selects Times Roman as basic font

\makeindex             % used for the subject index
\begin{document}

\title*{Domain Decomposition for the Closest Point Method}
\author{Ian May and Ronald D. Haynes and Steven J. Ruuth}
\institute{Ian May \at Simon Fraser University, 8888 University Dr, Burnaby, BC V5A 1S6, \email{mayianm@sfu.ca}
      \and Ronald Haynes \at Memorial University of Newfoundland, 230 Elizabeth Ave, St. Johns, NL A1C 5S7 \email{rhaynes@mun.ca}
      \and Steven Ruuth \at Simon Fraser University, 8888 University Dr, Burnaby, BC V5A 1S6, \email{sruuth@sfu.ca}}
\maketitle

\abstract*{The discretization of elliptic PDEs leads to large coupled systems of equations. Domain decomposition methods (DDMs) are one approach to the solution of these systems, and can split the problem in a way that allows for parallel computing. Herein, we extend two DDMs to elliptic PDEs posed intrinsic to surfaces as discretized by the Closest Point Method (CPM) \cite{SJR:CPM,CBM:ICPM}. We consider the positive Helmholtz equation
\begin{equation}
  \left(c-\Delta_\surf\right)u = f,
  \label{MAYDD25:eq:model}
\end{equation}
where $c\in\mathbb{R}^+$ is a constant and $\Delta_\surf$ is the Laplace-Beltrami operator associated with the surface $\surf\subset\mathbb{R}^d$. The evolution of diffusion equations by implicit time-stepping schemes and Laplace-Beltrami eigenvalue problems \cite{CBM:Eig} both give rise to equations of this form. The creation of efficient, parallel, solvers for this equation would ease the investigation of reaction-diffusion equations on surfaces \cite{CBM:RDonPC}, and speed up shape classification \cite{Reuter:ShapeDNA}, to name a couple applications.}

\section{Introduction}
\label{MAYDD25:intro}
The discretization of elliptic PDEs leads to large coupled systems of equations. Domain decomposition methods (DDMs) are one approach to the solution of these systems, and can split the problem in a way that allows for parallel computing. Herein, we extend two DDMs to elliptic PDEs posed intrinsic to surfaces as discretized by the Closest Point Method (CPM) \cite{SJR:CPM,CBM:ICPM}. We consider the positive Helmholtz equation
\begin{equation}
  \left(c-\Delta_\surf\right)u = f,
  \label{MAYDD25:eq:model}
\end{equation}
where $c\in\mathbb{R}^+$ is a constant and $\Delta_\surf$ is the Laplace-Beltrami operator associated with the surface $\surf\subset\mathbb{R}^d$. The evolution of diffusion equations by implicit time-stepping schemes and Laplace-Beltrami eigenvalue problems \cite{CBM:Eig} both give rise to equations of this form. The creation of efficient, parallel, solvers for this equation would ease the investigation of reaction-diffusion equations on surfaces \cite{CBM:RDonPC}, and speed up shape classification \cite{Reuter:ShapeDNA}, to name a couple applications.

Several methods exist for the discretization of surface intrinsic PDEs. The surface may be parametrized to allow the use of standard methods in the parameter space \cite{FloaterHormann:Para}. Unfortunately, many surfaces of interest do not have simple, or even known, parametrizations. Given a triangulation of the surface, a finite element discretization can be formed \cite{DziukElliot}. This approach leads to a sparse and symmetric system but is sensitive to the quality of the triangulation. Level set methods for surface PDEs \cite{Cheng:LSM} solve the problem in a higher dimensional embedding space over a narrow band containing the surface. The solution of model equation \eqref{MAYDD25:eq:model} by this method requires using gradient descent, as the approach was formulated only for parabolic problems. The CPM is also discretized over a narrow band in the embedding space, but has the advantage of using a direct discretization of equation \eqref{MAYDD25:eq:model}.

The solution of the linear system arising from the CPM discretization of the model equation \eqref{MAYDD25:eq:model} has relied primarily on direct methods, although a multigrid method was discussed in \cite{Chen:MG}. Herein we formulate restricted additive Schwarz (RAS) and optimized restricted additive Schwarz (ORAS) solvers compatible with the CPM to step towards efficient iterative solvers and to allow for parallelism. The optimized variant of the classical RAS solver uses Robin transmission conditions (TCs) to pass additional information between the subdomains \cite{Gand:OptPar,Cyr:OMSORAS}, and can accelerate convergence dramatically. This formulation is described in Sections \ref{MAYDD25:subd} and \ref{MAYDD25:trans} after reviewing the CPM in Section \ref{MAYDD25:cpm} and (O)RAS solvers in Section \ref{MAYDD25:oras}. Then, we discuss a PETSc \cite{petsc-user-ref,petsc-efficient} implementation and show some numerical examples in Section \ref{MAYDD25:res}. A more thorough exploration of these solvers, and an initial look at their use as preconditioners, can be found in May's thesis \cite{May:Thesis}.

\section{The closest point method}
\label{MAYDD25:cpm}
The CPM was introduced in \cite{SJR:CPM} as an embedding method for surface intrinsic PDEs. It allows the reuse of standard flat space discretizations of differential operators and provides a surface agnostic implementation. At the core of this method is the closest point mapping, $CP_\surf(x) = \underset{y\in\surf}{\argmin}|x-y|$ for $x\in\mathbb{R}^d$,
%\begin{eqnarray}
 %   CP_\surf:\mathbb{R}^d&\rightarrow&\surf \\
  %  x&\rightarrow&\underset{y\in\surf}{\argmin}|x-y| \nonumber
  %\label{MAYDD25:eq:cpmap}
%\end{eqnarray}
which identifies the closest point on the surface for (almost) any point in the embedding space. This mapping exists and is continuous in the subset of $\mathbb{R}^d$ consisting of all points within a distance $\kappa_\infty^{-1}$ of the surface, where $\kappa_\infty$ is an upper bound on the principal curvatures of the surface \cite{Chu:Vari}.

From this mapping, an extension operator $E$ can be defined that sends functions defined on the surface, $f:\surf\rightarrow\mathbb{R}$, to functions defined on the embedding space via composition with the closet point mapping, $Ef = f\circ CP_\surf$. The extended functions are constant in the surface normal direction and retain their original values on the surface. This extension operator can be used to define surface intrinsic differential operators from their flat space analogs \cite{SJR:CPM}.

Discretization typically requires a Cartesian grid on the embedding space within a narrow tube surrounding the surface. The extension operator can be defined by any suitable interpolation scheme, with tensor product barycentric Lagrangian interpolation \cite{Tref:Bary} being used here. As such, the computational tube must be wide enough to contain the interpolation stencil for any point on the surface. Using degree $p$ interpolation  and a grid spacing of $\Delta x$ requires that the tube contain all points within a distance of $\gamma=\Delta x(p+2)\sqrt{d}/2$ from the surface, thus limiting the acceptable grid spacings in relation to $\kappa_\infty$.
%\begin{equation}
%  \gamma = \frac{(p+2)\Delta x\sqrt{d}}{2}.
%  \label{MAYDD25:eq:tubeWidth}
%\end{equation}
%However, $\gamma\leq\kappa_\infty^{-1}$ must hold, thus constraining the grid spacing.

The grid points within the computational tube form the set of active nodes, %which we shall denote by 
$\Sigma_A$. For (\ref{MAYDD25:eq:model}), we need only discretize the regular Laplacian on $\mathbb{R}^d$. Here we consider the second order accurate centered difference approximation requiring $2d+1$ points. Around $\Sigma_A$ and lying outside the tube, a set of ghost nodes, $\Sigma_G$,  is formed from any incomplete differencing stencils. 
 With a total of $N_A$ active nodes and $N_G$ ghost nodes, we define the discrete Laplacian and extension operators, $\Delta^h:\mathbb{R}^{N_A+N_G}\rightarrow\mathbb{R}^{N_A}$ and $\mat{E}:\mathbb{R}^{N_A}\rightarrow\mathbb{R}^{N_A+N_G}$, 
%% \begin{eqnarray}
%%   \Delta^h:\mathbb{R}^{N_A+N_G}\rightarrow\mathbb{R}^{N_A}, \\
%%   \mat{E}:\mathbb{R}^{N_A}\rightarrow\mathbb{R}^{N_A+N_G},
%% \end{eqnarray}
where $\Delta^h$ applies the centered difference Laplacian over all active nodes, and $\mat{E}$ is the discretization of $E$. $\mat{E}$  extends data on the active nodes to both the active and ghost nodes, and has entries consisting of the interpolation weights for each node's closest point.

The Laplace-Beltrami operator can be directly discretized as $\Delta^h_{\surf,dir}=\Delta^h\mat{E}$, which was used successfully for parabolic equations with explicit time-stepping in \cite{SJR:CPM}. However, for implicit time-stepping \cite{CBM:ICPM} and eigenvalue problems \cite{CBM:Eig} a modified form is needed. In \cite{CBM:ICPM} it was recognized that there was a redundant interpolation being performed, and that its removal could stabilize the discretization. The stabilized form
%\begin{equation}
$  \Delta^h_\surf = -\frac{2d}{\Delta x^2}\mat{I} + \left(\frac{2d}{\Delta x^2}\mat{I} + \Delta^h\right)\mat{E},$
%  \label{MAYDD25:eq:stabLap}
%\end{equation}
will be used in the remainder of this work.

\section{(Optimized) Restricted additive Schwarz}
\label{MAYDD25:oras}
Both RAS and ORAS are overlapping DDMs, and can work on the same set of subdomains (given an additional overlap condition for ORAS \cite{Cyr:OMSORAS}). We define these solvers from the continuous point of view and subsequently discretize, rather than defining them purely algebraically.  This will ease the discussion of TCs within the context of the CPM later in Section \ref{MAYDD25:trans}.

First, the whole surface $\surf$ is decomposed into $N_S$ disjoint subdomains, $\widetilde{\surf}_j$, for $j=1,\ldots, N_S$. These disjoint subdomains are then grown to form overlapping subdomains $\surf_j$, whose boundaries are labelled depending on where they lie in the disjoint partitioning. Taking $\Gamma_{jk} = \partial\surf_j\cap\widetilde{\surf}_k$ gives $\partial\surf_j = \bigcup\limits_k\Gamma_{jk}$ and allows the definition of the local problems
\begin{equation}
  \begin{cases}
    \left(c-\Delta_\surf\right)u_j^{(n+1)} = f,\quad&{\rm in}~\surf_j, \\
    \mathcal{T}_{jk}u_j^{(n+1)} = \mathcal{T}_{jk}u^{(n)},\quad&{\rm on}~\Gamma_{jk},~k=1,\ldots, N_S,~k\neq j,
  \end{cases}
  \label{MAYDD25:eq:contSub}
\end{equation}
where $\mathcal{T}_{jk}$ are generally linear boundary operators defining the TCs. RAS is achieved by choosing $\mathcal{T}_{jk}$ as identity operators, corresponding to Dirichlet TCs, 
while ORAS uses Robin TCs, $\mathcal{T}_{jk} = \left(\pdiff{}{\hat{\mat{n}}_{jk}}+\alpha\right)$, where $\hat{\mat{n}}_{jk}$ is the outward pointing boundary normal on $\Gamma_{jk}$ and $\alpha\in\mathbb{R}^+$ is a constant weight on the Dirichlet contribution.

The subproblems in equation \eqref{MAYDD25:eq:contSub} are initialized with a guess for the global solution $u^{(0)}$ (defined at least over the boundaries $\Gamma_{jk},\forall j,k$), which is usually just taken as $u^{(0)}=0$. After all of the subproblems have been solved a new global solution is constructed with respect to the disjoint partitioning,
  $u^{(n+1)} = \sum\limits_j\restrict{u_j^{(n+1)}}{\widetilde{\surf}_j}$,
%\begin{equation}
%  u^{(n+1)} = \sum\limits_j\restrict{u_j^{(n+1)}}{\widetilde{\surf}_j},
%  \label{MAYDD25:eq:solCon}
%\end{equation}
where the use of the term \textit{restricted} indicates that the portion of the local solutions in the overlap regions are discarded. From this new approximation for the global solution the local problems may be solved again with new boundary data, and the process repeats until the global solution is satisfactory.

\section{Subdomain construction}
\label{MAYDD25:subd}
To solve problems arising from the CPM we first need to decompose the global set of active nodes $\Sigma_A$. (O)RAS solvers rely on both a disjoint partitioning of the active nodes and an induced overlapping partitioning.
%Here, we describe one approach generating this disjoint partitioning.
  %The growth of the disjoint partitions into overlapping subdomains, and the formation of some relevant geometric information, follows the same procedure regardless of how the disjoint partitioning is found.
%Throughout the following discussion we will rely on several sets of nodes.
Following the notation in Section \ref{MAYDD25:oras}, disjoint partitions will be denoted by $\widetilde{\Sigma}_j$, overlapping partitions by $\Sigma_j$, and the boundaries of the overlapping partitions by $\Lambda_j$.

\begin{figure}
  \vspace{-2em}
  \sidecaption
  \includegraphics[scale=0.45]{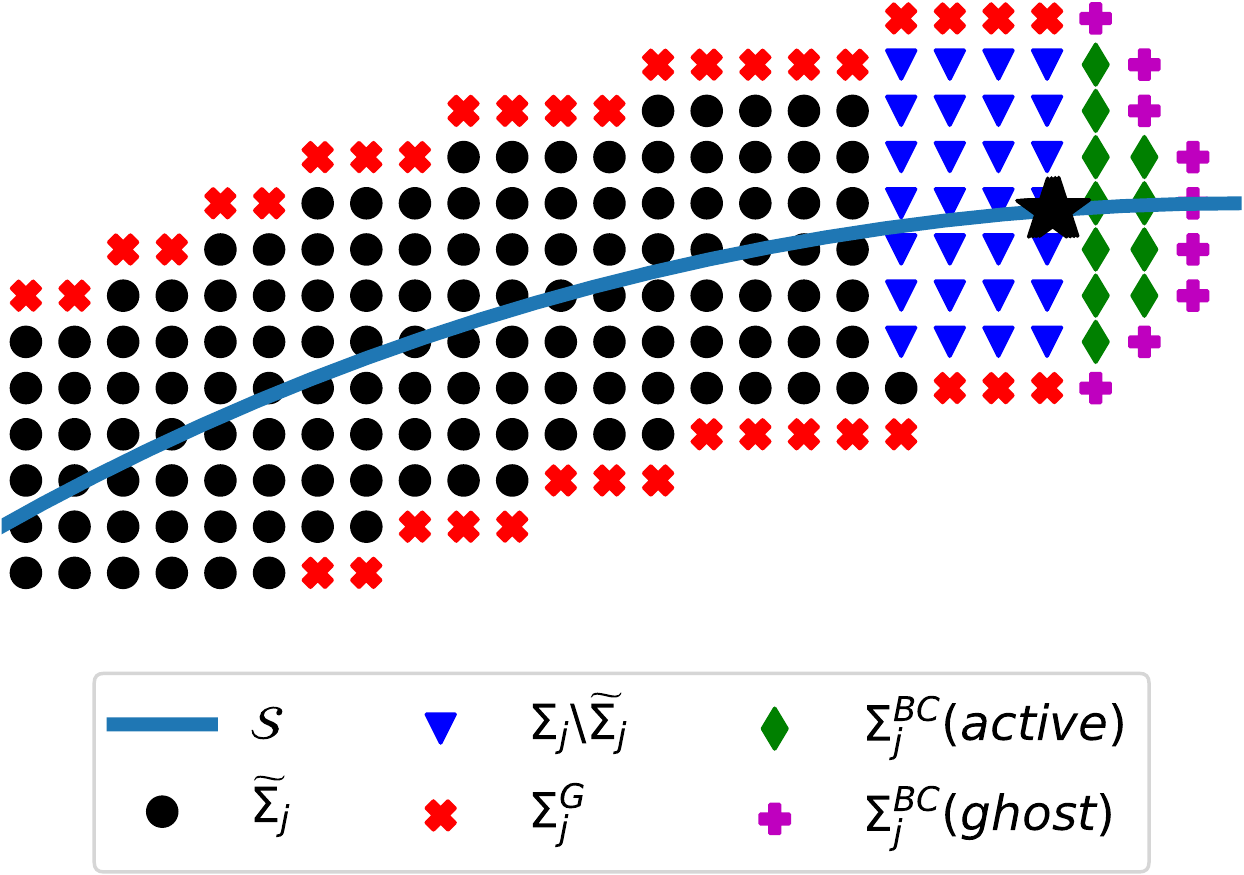}
  \caption{A portion of a subdomain from a circular problem with $N_S=8$ subdomains and $N_O=4$ layers of overlap is shown here. The nodes are marked according to their role in the subdomain as described in Section \ref{MAYDD25:subd}. The points belonging to the set $\Lambda_j$ are shown by (nearly coincident) black stars.}
  \label{MAYDD25:fig:scons}
  \vspace{-1em}
\end{figure}

To ensure the solvers work on a variety of surfaces, we seek an automated and surface agnostic partitioning scheme to generate the disjoint partitions. METIS \cite{metis} is a graph partitioner that is frequently used within the DD community to partition meshes \cite{DoleanNataf}. The stencils of $\Delta^h$ and $\mat{E}$ may be used to induce connectivity between the active nodes and define a graph. Here we only consider nearest neighbor coupling through the stencil for $\Delta^h$. Fig. \ref{MAYDD25:fig:scons} shows a portion of one such disjoint partition, in black circles,  for a circular surface. 

With $\widetilde{\Sigma}_j$ obtained from METIS, overlapping subdomains $\Sigma_j$ can be formed. This construction proceeds in the following steps:
\begin{enumerate}
\item All nodes in $\widetilde{\Sigma}_j$ are added to $\Sigma_j$.
\item $N_O$ layers of overlap nodes are added around $\Sigma_j$. Layers are added one at a time from globally active nodes neighboring $\Sigma_j$.
\item A subset of the ghost nodes, $\Sigma_G$, are placed in $\Sigma_j^G$ which consists of nodes that neighbor a member of $\Sigma_j$.
\item The shapes of the disjoint and overlapping subdomains are not known in advance. The boundary $\partial\surf_j$ is approximated discretely by the closest points of the final layer of overlap nodes, and held in the set $\Lambda_j$.
\item Nodes needed to complete stencils from the ambient Laplacian or extension operator, including extension from the points $x_i\in\Lambda_j$, are placed in the set $\Sigma_j^{BC}$.
\item For ORAS a layer of ghost nodes around $\Sigma_j^{BC}$ are also placed in $\Sigma_j^{BC}$.
\end{enumerate}

The active nodes in the $j^{\rm th}$ subdomain consist of $\Sigma_j$ and the active portion of $\Sigma_j^{BC}$. $\Sigma_j^{BC}$ is kept separate as that is where the TCs in Section \ref{MAYDD25:trans} are defined. Each of these sets are shown in Fig. \ref{MAYDD25:fig:scons}, which shows a portion of one subdomain on a circle in the vicinity of the points in $\Lambda_j$ at one of its boundaries.

The Robin TCs, to be defined in Section \ref{MAYDD25:trans}, need some final information about the subdomain. Every node in $\Sigma_j^{BC}$ is identified with the point in $\Lambda_j$ that is closest to it. This identification will be used to override the global closest point function in the following section. For each point in $\Lambda_j$ we also need to know the direction that is simultaneously orthogonal to the boundary and the surface normal direction. We call this the conormal direction. It is in this direction that the Neumann component of the Robin condition will be enforced. However, the discrete nature of $\Lambda_j$ makes this construction difficult. Instead we define the conormal vectors from the point of view of the boundary nodes. Take $x_i\in\Sigma_j^{BC}$ as a node whose associated conormal direction, $\hat{q}_i$, is sought. Let $y_i$ be its closest point in $\Lambda_j$, and $\hat{n}_i$ be the unit surface normal there. Connecting the boundary location to the boundary node via $d_i = x_i - y_i$, we obtain a usable approximation to the conormal by computing the component of $d_i$ that is orthogonal to $\hat{n}_i$ and normalizing, i.e., $\hat{q}_i = (d_i - \left(d_i\cdot\hat{n}_i\right)\hat{n}_i)/(\left|d_i - \left(d_i\cdot\hat{n}_i\right)\hat{n}_i\right|).$ In the unlikely event that $d_i$ lies perfectly in the surface normal direction, we set $\hat{q}_i=0$ which recovers the natural boundary condition on the computational tube as discussed at the end of Section \ref{MAYDD25:trans}.

\section{Transmission conditions}
\label{MAYDD25:trans}
Boundary conditions in the CPM are imposed by modifying the extension operator over the nodes $\Sigma_j^{BC}$ beyond the surface boundary \cite{CBM:Eig}. %The nodes in $\Sigma_j^{BC}$ are precisely the nodes needing this modification, and this is why they were kept separate in the above construction. 
As such, the local operators will take the form
\begin{equation}
  \mat{A}_j = \left(c+\frac{2d}{h^2}\right) - \left(\frac{2d}{h^2}+\Delta_j^h\right)\begin{bmatrix}\mat{E}_j \\ \mat{T}_j\end{bmatrix},
  \label{MAYDD25:eq:localOp}
\end{equation}
where $\mat{E}_j$ is the extension operator for the nodes in $\Sigma_j$ as inherited from the global operator and $\mat{T}_j$ is the modified extension operator for the nodes in $\Sigma_j^{BC}$. When solving for the local correction to the solution the right hand side of the local problem, $\mat{A}_jv_j=r_j$, will be the restriction of the residual to $\Sigma_j$. The final rows of the right hand side, those lying over $\Sigma_j^{BC}$, become zeros corresponding to the homogenous TCs. 
%also need to be changed in accordance with the TC. These conditions become homogeneous, and the right hand side is simply zero.

% Dirfo --------------------------------
Homogeneous Dirichlet TCs can be enforced to first order accuracy by extending zeros over all of $\Sigma_j^{BC}$. With the right hand side already set to zero there, the modified extension reduces to the identity mapping, $\mat{T}_j = \begin{bmatrix} {\bf 0} & ~ & {\bf I}\end{bmatrix}$, with the zero matrix padding the columns corresponding to the interior nodes.

% Robfo -------------------------------
We discretize the Robin condition
\begin{equation}
  \restrict{\pdiff{u}{\hat{q}_i}}{CP_{\surf_j}(x_i)} + \alpha u\left(CP_{\surf_j}(x_i)\right) = 0,
  \label{MAYDD25:eq:contRob}
\end{equation}
using a forward difference in the $\hat{q}_i$ direction for each node in $\Sigma_j^{BC}$ and the first order accurate Dirichlet condition from above. Taking the partial derivative $\pdiff{u}{d_i}$, and applying the change of variables $d_i = \hat{q}_i + \hat{n}_i$, allows one to write the Neumann term in equation \eqref{MAYDD25:eq:contRob} in terms of the displacement vector $d_i$ from Section \ref{MAYDD25:subd}. Assuming for the moment that $d_i$ and $\hat{q}_i$ are not perpendicular, the derivative in the conormal direction can be approximated by
%\begin{equation}
$  \restrict{\pdiff{u}{\hat{q}_i}}{CP_{\surf_j}(x_i)} \approx \frac{u(x_i) - u\left(CP_{\surf_j}(x_i)\right)}{d_i\cdot\hat{q}_i}$
%  \label{MAYDD25:eq:conormDiff}
%\end{equation}
where $CP_{\surf_j}$ denotes the modified closest point function identifying points in $\Sigma_j^{BC}$ with points in $\Lambda_j$. 
Combining this with (\ref{MAYDD25:eq:contRob}), and applying the identity extension for the Dirichlet component, $u(x_i)=u(CP_{\surf_j}(x_i))$, we find that $\mat{T}_j$ must enforce the extension
%
%After adding the identity extension for the Dirichlet component, $\mat{T}_j$ must enforce the extension
%\begin{equation}
$  u(x_i) = \frac{u\left(CP_{\surf_j}(x_i)\right)}{1 + \alpha d_i\cdot\hat{q}_i},$
%  \label{MAYDD25:eq:robfo}
%\end{equation}
with $u\left(CP_{\surf_j}(x_i)\right)$ replaced by the same interpolation used in the global scheme discussed in Section \ref{MAYDD25:cpm}.

As $d_i$ approaches the surface normal direction, $d_i\cdot\hat{q}_i$ will tend to zero. In this event, the extension reduces to $u(x_i) = u\left(CP_{\surf_j}(x_i)\right)$, which is just the standard extension corresponding to the interior. Fortuitously, this case arises when the point $x_i$ lies adjacent to the interior points where this condition would be applied anyway, and in our experience this ensures that the method remains robust.

\section{Results}
\label{MAYDD25:res}
The solvers described in the previous sections were implemented in C\raisebox{0.5ex}{\small\textbf{++}}, with PETSc \cite{petsc-user-ref,petsc-efficient} providing the linear algebra data structures and MPI parallelization, and Umfpack \cite{Davis:Umfpack} providing the local solutions. Here we focus on evaluating the solver, though in practice one should accelerate the solver with a Krylov method. The (O)RAS solver was placed into a PETSc \texttt{PCSHELL} preconditioner, allowing it to be embedded in any of their Krylov methods, and we have found coupling with GMRES to be a favorable pair.

Equation \eqref{MAYDD25:eq:model} was solved over the Stanford Bunny \cite{bunny}, which has been scaled to be two units tall. The original triangulation has not been modified in any way beyond this scaling. This surface has several holes and is complicated enough to stress the solvers, making it a good test case. Our chosen grid spacing was $\Delta x = 1/120$, which paired with tri-quadratic interpolation gives $N_A=947,964$ active nodes in the global problem. The origin was placed at the center of the bounding box containing the bunny and the right hand side $f=\phi(\pi-\phi)\sin(3\phi)(\sin\theta+\cos(10\theta))/2$ was used after extending it to be constant along the surface normals.

\begin{table}
  \centering
  \begin{tabular}{@{}rccc|rccc|rccc@{}}
    \multicolumn{4}{c|}{$N_O=4$, $\alpha=16$} & \multicolumn{4}{c|}{$N_S=64$, $\alpha=16$} & \multicolumn{4}{c}{$N_S=64$, $N_O=4$} \\ \midrule
    $N_S$    & $64$     & $96$    & $128$     & $N_O$   & $4$       & $6$      & $8$       & $~\alpha~$ & $16$  & $32$  & $64$     \\
    RAS      & $992$    & $1237$  & $1533$    &         & $992$     & $747$    & $610$     &            & $992$ & $992$ & $992$    \\
    ORAS     & $526$    & $672$   & $833$     &         & $526$     & $418$    & $393$     &            & $526$ & $707$ & $868$
  \end{tabular}
  \caption{Here the iterations to convergence of the ORAS solver are gathered for  various parameters. Convergence was declared when the $2-$norm of the residual was reduced by a factor of $10^6$.}
  \label{MAYDD25:tab:full}
\end{table}

\begin{figure}
  \centering
  \includegraphics[width=0.9\textwidth]{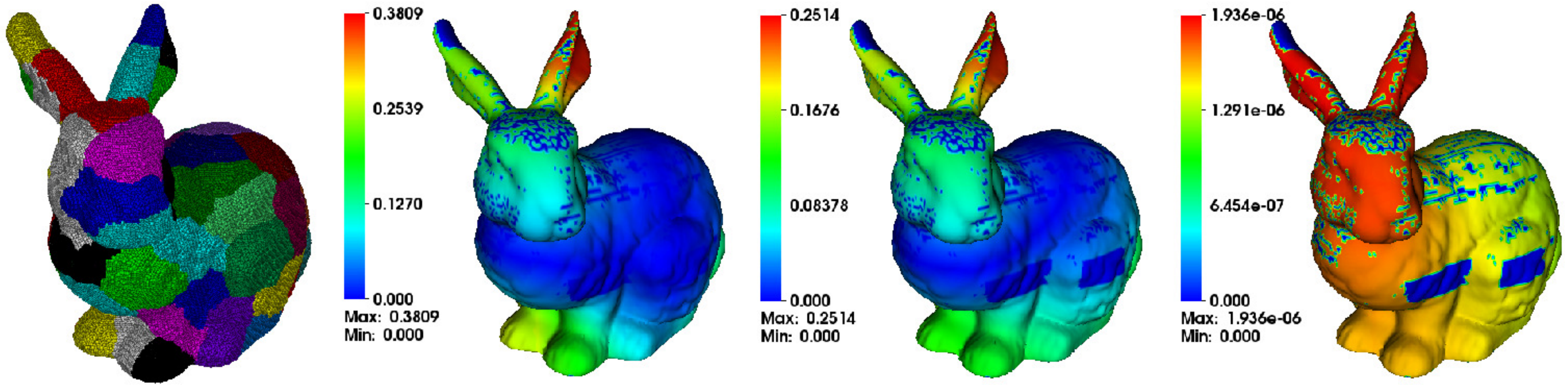}
  \caption{The Stanford Bunny test problem solved with ORAS using $N_S=64$, $N_O=4$, $\alpha=16$. The first panel shows the disjoint partitioning from METIS. The second, third, and fourth panels show the error in the solution after the $1^{\rm st}$, $10^{\rm th}$, and $500^{\rm th}$ iterations compared to the converged solution.}
  \label{MAYDD25:fig:bunny}
  \vspace{-3em}
\end{figure}

Table \ref{MAYDD25:tab:full} shows the effects of subdomain count $N_S$, overlap width $N_O$, and Robin parameter $\alpha$. For comparison, GMRES preconditioned with the standard block-Jacobi method with $64$, $96$, and $128$ blocks requires more than $10000$ iterations. The solvers display the expected behavior with the iteration count increasing for larger subdomain counts and decreasing with larger overlap widths. ORAS consistently requires fewer iterations than RAS, though the final sub-table shows the dependence of this performance on the appropriate choice of Robin weight $\alpha$. The partitioning, the initial error, and the error in the approximate solution after $10$ and $550$ iterations are visible in Fig. \ref{MAYDD25:fig:bunny} for one run of the solver.

Choosing an optimal value for $\alpha$ is non-trivial as it depends on the value of $c$, the mesh width, and the geometry. Additionally, the presence of cross points in decomposition, where more than two subdomains meet, complicate the matter. From the planar case, it is known that $\alpha\sim\mathcal{O}\left(\Delta x^{-1/2}\right)$, but determining precise values \textit{a priori} is limited to simple splittings \cite{Gand:OptPar,Gand:XPT}. An upcoming work from the same authors explores this in much greater detail.

\section{Conclusion}
\label{MAYDD25:conc}
Restricted additive Schwarz and optimized restricted additive Schwarz solvers were formulated for the closest point method applied to \eqref{MAYDD25:eq:model}. These solvers provide a solution mechanism for larger problem sizes and will allow users of the CPM to leverage large scale parallelism. Table \ref{MAYDD25:tab:full} shows the dramatic reduction in iteration count when Robin TCs are used. These solvers were more completely evaluated in \cite{May:Thesis}, which includes an exploration of their utility as preconditioners. The optimized conditions come at the cost of some additional complexity in the implementation, and even the standard RAS solver brings parallel capabilities to the user. Interesting extensions to this work include multiplicative methods, non-overlapping Robin schemes, two-level solvers, and inclusion of advective terms in the model equation.

{\em Acknowledgements} The authors gratefully acknowledge the financial support of NSERC Canada (RGPIN 2016-04361 and RGPIN 2018-04881), and the preliminary work of Nathan King that inspired this project.

% Bibliography
\bibliographystyle{spmpsci}
\bibliography{references}

\end{document}